\documentclass[12pt,a4paper,twoside]{article}
\usepackage{amssymb,amsfonts,amsmath,amsthm,euscript}
\usepackage[english]{babel}
\usepackage[latin1]{inputenc}
\usepackage{dsfont}
\usepackage{color}
\usepackage{placeins}
\usepackage{mathrsfs}
\usepackage{geometry}
\usepackage{graphicx}
\usepackage{subfigure}
\usepackage{mathtools}
\usepackage{multirow}
\usepackage[bf,small]{caption}
\usepackage{sectsty}
\sectionfont{\large}
\begin{document}


\theoremstyle{plain}
\newtheorem{teo}{Theorem}[section]
\newtheorem{corol}{Corollary}[section]
\newtheorem{prop}{Proposition}[section]
\newtheorem{lema}{Lemma}[section]
\theoremstyle{definition}
\newtheorem{exem}{Example}[section]
\newtheorem{defin}{Definition}[section]
\newtheorem{obs}{Remark}[section]

\renewcommand{\theequation}{\thesection.\arabic{equation}}
\renewcommand{\thetable}{\thesection.\arabic{table}}
\renewcommand{\thefigure}{\thesection.\arabic{figure}}
\renewcommand{\arraystretch}{1.2}
\long\def\symbolfootnote[#1]#2{\begingroup%
\def\thefootnote{\fnsymbol{footnote}}\footnote[#1]{#2}\endgroup}
\numberwithin{table}{section}
\numberwithin{figure}{section}


%
%

\newcommand{\cov}{\mathrm{Cov}}
\newcommand{\corr}{\mathrm{Corr}}
\newcommand{\C}{\mathbb{C}}
\newcommand{\E}{\mathbb {E}}
\newcommand{\F}[1]{\mathcal{F}_{#1}}
\newcommand{\iten}[1]{\vspace{.2cm}\noindent{\bf #1}}
\newcommand{\K}{\mathcal{K}}
\renewcommand{\L}{\ell}
\newcommand{\B}{\mathcal{B}}
\newcommand{\N}{\mathbb {N}}
\newcommand{\nn}[1]{\mbox{\boldmath{$#1$}}}
\newcommand{\noise}{$\{Z_t\}_{t \in \mathbb{Z}}{}$}
\newcommand{\prob}{\mathbb{P}}
\newcommand{\pe}{$\{X_t\}_{t \in \mathbb{Z}}{}$}
\renewcommand{\proof}{\noindent\textbf{Proof: }}
\renewcommand{\qed}{\begin{flushright}\vspace{-.5cm}$\Box$\end{flushright}}
\newcommand{\R}{\mathbb {R}}
\newcommand{\stg}{$\left\{X_{t}\right\}_{t=1}^{n}{}$}
\newcommand{\var}{\mathrm{Var}}
\newcommand{\Z}{\mathbb {Z}}
\newcommand{\mat}[1]{\mbox{\boldmath{$#1$}}}
\newcommand{\fim}{\hfill $\Box$}
\newcommand{\stk}[1]{\stackrel{#1}{\longrightarrow}}
\newcommand{\s}{\hphantom{x}}

\newcommand{\I}{\mathds{1}}
\thispagestyle{empty}

\vskip2cm
{\centering
\Large{\bf Fractionally Integrated Moving Average Stable Processes With Long-Range Dependence}	

\vspace{1.0cm}
\large{\bf{G.L. Feltes and S.R.C. Lopes \symbolfootnote[3]{Corresponding author. E-mail: silviarc.lopes@gmail.com}}}

\vspace{0.2cm}
Mathematics and Statistics Institute \\
Federal University of Rio Grande do Sul\\
Porto Alegre - RS - Brazil \\
}

\vspace{0.6cm}

\centerline{\today}

\vspace{0.4cm}


\begin{abstract}

\noindent Long memory processes driven by L\'evy noise with finite second-order moments have been well studied in the literature. They form a very rich class of processes presenting an autocovariance function that decays like a power function. Here, we study a class of L\'evy processes whose second-order moments are infinite, the so-called $\alpha$-stable processes. Based on Samorodnitsky and Taqqu (2000), we construct an isometry that allows us to define stochastic integrals concerning the linear fractional stable motion using Riemann-Liouville fractional integrals. With this construction, an integration by parts formula follows naturally. We then present a family of stationary $S\alpha S$ processes with the property of long-range dependence, using a generalized measure to investigate its dependence structure. At the end, the law of large number's result for a time's sample of the process is shown as an application of the isometry and integration by parts formula.

\vspace{0.2cm}
\vspace{.2cm}
\end{abstract}

\section{Introduction}\label{sec_introduction}
\renewcommand{\theequation}{\thesection.\arabic{equation}}
\setcounter{equation}{0}

In this work, we are interested in continuous-time fractionally integrated stable processes and its property of long-range dependence (LRD) or long memory. There are many interesting applications of such processes throughout the literature (see Samorodnitsky, 2006 for an intuitive and historical introduction). This property presents a phenomenon in which some proper measure of the process' dependence structure - usually the autocovariance function when well defined - behaves like a power function, having a slow decay. Slow in the sense that the corresponding sum will diverge (power function with exponent in $(-1,0)$). Linear fractional stable motion (LFSM) is a generalization of the well-known fractional Brownian motion (fBm) to the $\alpha$-stable case. Accordingly to Samorodnitsky and Taqqu (1994), it has the following moving average representation. For $0 < \alpha \le 2$, $0 < H < 1$, $H \neq \frac{1}{\alpha}$ and $a,b \in \R$ such that $|a| + |b| > 0$,

\begin{equation} \label{I_eq_1}
M_{\alpha,H}(t) = \int_{-\infty}^{\infty} a\left[(t - x)_+^{H - \frac{1}{\alpha}} - (-x)_+^{H - \frac{1}{\alpha}}\right] + b\left[(t - x)_-^{H - \frac{1}{\alpha}} - (-x)_-^{H - \frac{1}{\alpha}}\right] dL_\alpha(x) \; ,
\end{equation}

\noindent for $t \in \R$, where $u_+ = \max\{u,0\}$ and $u_- = \max\{-u,0\}$ are the positive and negative part of $u$, respectively, and $\{L_\alpha(t)\}_{t \in \R}$ is a $S\alpha S$ process with Lebesgue control measure. The process $\{M_{\alpha,H}(t)\}_{t \in \R}$ is self-similar with Hurst parameter $H$ and it has stationary increments. Defining the increment process $Y_{\alpha,H}(t) := M_{\alpha,H}(t + 1) - M_{\alpha,H}(t) \; , \; t \in \R$ which is also stationary, consider as a measure of its dependence structure the difference between the joint characteristic function of the pair $(Y_{\alpha,H}(t), Y_{\alpha,H}(0))$ and the product of the marginal characteristic functions of $Y_{\alpha,H}(t)$ and $Y_{\alpha,H}(0)$, given by $$r(t) := \E\exp\{i(\theta_1Y_d(t) + \theta_2Y_d(0))\} - \E\exp\{i\theta_1Y_d(t)\}\E\exp\{i\theta_2Y_d(0)\} \; .$$ Astrauskas et al. (1991) show that, for a positive constant $K$, as $t \to \infty$,
\begin{equation*}
	\left|r(t)\right|	\sim K|t|^{\alpha H - \alpha} \; ,
\end{equation*}
\noindent and when $1 < \alpha < 2$ and $\frac{1}{\alpha} < H < 1$ we have $-1 < \alpha H - \alpha < 0$. Therefore, one can say that the increment process of the LFSM has the property of long-range dependence. This fact suggests that we may be able to define a continuous-time $\alpha$-stable process with long memory if we succeed in building stochastic integrals with respect to the LFSM.

\noindent This paper is organized as follows. Section 2 starts with some basic facts and introductory notation. Using Riemann-Liouville fractional integrals we build in Section 3 an isometry which allows us to define the stochastic integral with respect to an LFSM and prove an interesting formula for representing the fractionally integrated $\alpha$-stable processes. In Section 4 we define a long memory fractionally integrated moving average $\alpha$-stable process and relate it to the LFSM using previous results. As an application of the isometry and integration by parts formula, we show a law of large number's results at the end of this section. Section 5 concludes this work.

\section{Definitions and Preliminary Results}
\renewcommand{\theequation}{\thesection.\arabic{equation}}
\setcounter{equation}{0}

\noindent Let $\{L_\alpha(t)\}_{t \in \R}$ be a $S\alpha S$ process with index of stability $0 < \alpha \le 2$ and Lebesgue control measure $dm(x) = dx$. Accordingly to Samorodnitsky and Taqqu (1994), for a measurable function $f : \R \to \R$ such that $f \in L^\alpha(\R)$ the stable stochastic integral
\begin{equation}
	I(f) = \int_{\R} f(x) dL_\alpha(x)
\end{equation}
\noindent is well defined and has the following properties.

\begin{prop} \label{DPR_prop_1}
	The distribution of $I(f)$ is given by
	\begin{equation*}
		I(f) \sim S_\alpha(\sigma_f,0,0) \; ,
	\end{equation*}
	where $\sim$ denotes equality in distribution and
	\begin{equation}\label{norm}
		\sigma_f = \left(\int_{\R} |f(x)|^\alpha dx\right)^\frac{1}{\alpha} = ||f||_\alpha,
	\end{equation}
	which means that
	\begin{equation*}
		\E\left[\exp\left\{i\theta I(f)\right\} \right] = \exp\left\{-|\theta|^\alpha \int_{\R} |f(x)|^\alpha dx \right\}.
	\end{equation*}

\end{prop}
\noindent \textbf{Proof.}
	See proposition 3.4.1 in Samorodnitsky and Taqqu (1994).
\vspace{-0.2cm}
\qed

\noindent In order to define the fractionally integrated stable processes, we need to introduce the fractional integrals and derivatives of Riemann-Liouville (see Samko et al, 1993).

\begin{defin}
	Let $d \in (0,1)$. For $f \in L^p(\mathbb{R})$, where $1 \leq p < \frac{1}{d}$, we define
	\begin{equation}\label{rlRL}
	(I_-^df)(x) = \frac{1}{\Gamma(d)} \int_{x}^{\infty} f(t)(t-x)^{d-1} dt
\ \ \mbox{and} \ \
	(I_+^df)(x) = \frac{1}{\Gamma(d)} \int_{-\infty}^{x} f(t)(x-t)^{d-1} dt \; ,
	\end{equation}
respectively, the \textit{right and left-sided Riemann-Liouville fractional integrals of order $d$}, where $\Gamma(x) = \int_{0}^{\infty} t^{x-1}e^{-t} dt$ is the Gamma function.
\end{defin}

\noindent On the other hand, fractional differentiation was introduced as the inverse of fractional integration. Consider $d \in (0,1)$, $1 \leq p < \frac{1}{d}$ and denote by $I_\pm^d(L^p)$ the class of functions $\phi \in L^p$ that can be represented as a fractional integral $I_\pm^d f$ of some function $f \in L^p$.

\begin{defin}
	Let $\phi \in I_\pm^d(L^p)$ be of the form $\phi = I_\pm^df$, for a function $f \in L^p$. Then $f$ is precisely the \textit{right or left-sided Riemann-Liouville fractional derivative of order $d$}, respectively, given by
	\begin{equation}
	(\mathcal{D}_-^d\phi)(x) = -\frac{1}{\Gamma(1-d)} \frac{d}{dx}\int_{x}^{\infty} \phi(t)(t-x)^{-d} dt 	
\end{equation}
	and
	\begin{equation}
	(\mathcal{D}_+^d\phi)(x) = \frac{1}{\Gamma(1-d)} \frac{d}{dx}\int_{-\infty}^{x} \phi(t)(x-t)^{-d} dt,
	\end{equation}
	where $\frac{d}{dx}$ is the usual operation of differentiation with respect to $x$ and the convergence of the integrals at the singularity $t = x$ holds pointwise with probability 1, if $p=1$ and in the $L^p$ sense, if $p>1$.
\end{defin}

\begin{obs} \label{DPR_remark_1}
	In fractional processes as the LFSM defined in (\ref{I_eq_1}), for $0 < \alpha \le 2$, one has an interesting relationship between the fractional parameter $d$ and the Hurst parameter $H$ of self-similarity which is given by
	\begin{equation}
		H = d + \frac{1}{\alpha}.
	\end{equation}
	As we are more interested in the fractional integration point of view, we will often mention the parameter $d = H - \frac{1}{\alpha}$.
\end{obs}

\section{Main Results I}
\renewcommand{\theequation}{\thesection.\arabic{equation}}
\setcounter{equation}{0}
This section is based on Pipiras and Taqqu (2000) and Marquardt (2006) in which we define stochastic integrals with respect to the linear fractional stable motion (LFSM) and prove a formula that can be seen as integration by parts' formula, analogous to proposition 5.5 on Marquardt (2006) for the fractional Brownian motion. We restrict ourselves to the case $1 < \alpha < 2$ since we are interested in processes with the property of long-range dependence. While the cited authors deal with L\'evy processes with a finite second moment, we handle the infinite second-order moments of $S\alpha S$ processes by using their finite lower-order moments in a suitable space.

\subsection{Fractionally Integrated Processes and the $\bf{{\emph{B}}_{{\mat{\footnotesize{\alpha}}, \bf{\emph{p}}}}}$ Space}

\indent Let $\{X(t)\}_{t \in \R}$ be a $S\alpha S$ process, with index of stability $1 < \alpha < 2$, given in its integral representation, accordingly to the family $\{f_t\}_{t \in \R}$ of functions such that, for every $t \in \R$, $f_t: \R \to \R$, $f_t \in L^\alpha(\R)$ and $\{L_\alpha(t)\}_{t \in \R}$ is a $S\alpha S$ process with Lebesgue control measure $dm(x) = dx$. Then, according to Proposition \ref{DPR_prop_1}, we know that
\begin{equation} \label{MRI_eq_1}
	X(t) = \int_{\R} f_t(x) dL_\alpha(x) \sim S_\alpha(\sigma_t, 0, 0)
\end{equation}
defined and its distribution is totally determined by
\begin{equation}
	\sigma_t = \left(\int_{\R} |f_t(x)|^\alpha dx\right)^\frac{1}{\alpha} = ||f_t||_\alpha \; .
\end{equation}

\noindent It would be important to have an isometry property for the $S\alpha S$ processes as the Itô's isometry for processes with finite second-order moment. Since the $\alpha$-stable distributions only have finite moments of order $p < \alpha$, we shall consider the following known result.

\begin{prop} \label{MRI_prop_1}
	Let $\{X(t)\}_{t \in \R}$ be the process defined in (\ref{MRI_eq_1}) and let $1 \leq p < \alpha$. Then
	\begin{equation}
	\E[|X(t)|^p] = \E[|L_\alpha(1)|^p] \; ||f_t||_\alpha^p.
	\end{equation}
\end{prop}

\noindent \textbf{Proof.}
	From expression (\ref{MRI_eq_1}) and since for any positive $a$, $S_\alpha(a\sigma,0,0) \sim aS_\alpha(\sigma,0,0)$, where $\sim$ denotes equality in distribution, we have
	\begin{equation*}
		X(t) \sim S_\alpha(\sigma_t,0,0) \sim \sigma_tS_\alpha(1,0,0) \sim \sigma_tL_\alpha(1).
	\end{equation*}
	Therefore
	\begin{equation*}
		\E[|X(t)|^p] = \E[|\sigma_tL_\alpha(1)|^p] = \E[|L_\alpha(1)|^p] \; ||f_t||_\alpha^p.
	\end{equation*}
\qed
\noindent Now, for measurable functions  $g: \R \to \R$, $g \in L^1(\R)$ and $d \in (0, 1 - \frac{1}{\alpha})$, consider the right-sided Riemann-Liouville fractional integral of order $d$, $I_-^dg = \frac{1}{\Gamma(d)} \int_{x}^{\infty} g(t)(t-x)^{d-1} dt$ and denote by $\tilde{B}_\alpha$ the set
\begin{equation}
	\tilde{B}_\alpha = \bigg\{g \in L^1(\R) : \int_{-\infty}^{\infty} |(I_-^dg)(x)|^\alpha dx < \infty  \bigg\} \; .
\end{equation}
The set $\tilde{B}_\alpha$ can be seen as the set of kernel functions $g \in L^1(\R)$ for which the stable stochastic integral of its respective fractional integral $\int_{\R} (I_-^dg)(x) dL_\alpha(x)$ is well defined.
For the next proposition, we will need the following lemma (Leach, 1956).

\begin{lema}
	Let $1 < p < \infty$ and $1 < q < \infty$ such that $\frac{1}{p} + \frac{1}{q} = 1$. For a measurable function $f:\R \to \R$, $f \in L^p(\R)$ if and only if there exists a constant $C>0$ such that
	\begin{equation}
	\int_{\R} |f(x)h(x)| dx \leq C||h||_q \; ,
	\end{equation}
	for every function $h \in L^q(\R)$.
\end{lema}
\noindent Note that since $\alpha$ lies in the interval $(1,2)$, $L^\alpha(\R)$ is a Banach space. Using the previous lemma, we have the following.
\begin{prop} \label{MRI_prop_2}
	If $g \in L^1(\R)\cap L^\alpha(\R)$, then $g \in \tilde{B}_\alpha$.
\end{prop}
\noindent \textbf{Proof.}
	For every $h \in L^\beta(\R)$, such that $\frac{1}{\alpha} + \frac{1}{\beta} = 1$, we have
	\begin{align}\label{MRI_eq_2}
\nonumber
		\int_{\R} \left|h(x)(I_-^dg)(x)\right| dx &\le \int_{-\infty}^{\infty} \int_{x}^{\infty} \left|h(x)g(t)\frac{(t-x)^{d-1}}{\Gamma(d)}\right| \;dt \;dx \\
			&= \frac{1}{\Gamma(d)}\int_{-\infty}^{\infty} \int_{0}^{\infty} \left|h(x)g(s+x)s^{d-1}\right| \;ds\;dx.
	\end{align}
	Rewriting the last expression, we get
	\begin{equation}\label{MRI_eq_3}
		\int_{-\infty}^{\infty} \int_{0}^{\infty} \left|h(x)g(s+x)s^{d-1}\right| \;ds\;dx = I_1 + I_2,
	\end{equation}
where
\begin{equation*}\label{MRI_eq_4}
I_1 = \int_{-\infty}^{\infty}\int_1^{\infty}|h(x)g(s+x)s^{d-1}|\, ds \, dx \ \ \mbox{and} \ \
I_2 = \int_{-\infty}^{\infty}\int_0^1|h(x)g(s+x)s^{d-1}|\, ds \, dx.
\end{equation*}
\noindent By using Fubini's theorem and Holder's inequality we get
	\begin{equation} \label{MRI_eq_6}
		I_2 = \int_{0}^{1} s^{d-1} \int_{-\infty}^{\infty} |h(x)g(s+x)| dx\;ds
\le \int_{0}^{1} s^{d-1} ||h||_\beta ||g||_\alpha ds = \frac{1}{d}||h||_\beta ||g||_\alpha.
	\end{equation}
	For $I_1$ in (\ref{MRI_eq_4}), set $t = s+x$ and apply Holder's inequality to get
	\begin{align}\label{MRI_eq_7}
\nonumber
		I_1 & = \int_{-\infty}^{\infty} |g(t)| \int_{1}^{\infty} |h(t-s)s^{d-1}| ds\;dt \le \int_{-\infty}^{\infty} |g(t)|||h||_\beta \left(\int_{1}^{\infty} s^{\alpha(d-1)} ds\right)^{\frac{1}{\alpha}} dt  \\
& = \int_{-\infty}^{\infty} |g(t)|||h||_\beta \frac{1}{(\alpha(1-d) - 1)^\frac{1}{\alpha}} \;dt = \frac{1}{(\alpha(1-d) - 1)^\frac{1}{\alpha}} ||g||_1 ||h||_\beta \; .
	\end{align}
	Combine (\ref{MRI_eq_2}), (\ref{MRI_eq_3}), (\ref{MRI_eq_6}) and (\ref{MRI_eq_7}) to finally get
	\begin{equation}\label{MRI_eq_8}
	\int_{\R} \left|h(x)(I_-^dg)(x)\right| dx \le \left[\frac{1}{\Gamma(d)}\left(\frac{||g||_\alpha}{d} + \frac{||g||_1}{(\alpha(1-d) - 1)^\frac{1}{\alpha}}\right)\right] ||h||_\beta,
	\end{equation}
	where
		$$C = \frac{1}{\Gamma(d)}\left(\frac{||g||_\alpha}{d} + \frac{||g||_1}{(\alpha(1-d) - 1)^\frac{1}{\alpha}}\right) > 0.$$
\qed

\begin{obs} \label{MRI_remark_1}
	As well as for the right-sided Riemann-Liouville fractional integral $I_-^dg$, the same is true for the left-sided Riemann-Liouville fractional integral $I_+^dg$, which means that if $g \in L^1(\R)\cap L^\alpha(\R)$, then $\int_{\R} |(I_+^dg)(x)|^\alpha dx < \infty$ and its stable stochastic integral $\int_{\R} (I_+^dg)(x) dL_\alpha(x)$ is well defined. The proof is completely analogous and will be omitted.
\end{obs}

\noindent As we have seen, $L^1(\R) \cap L^\alpha(\R) \subset \tilde{B}_\alpha$, so we would like to define a norm which could allow us to work with isometries for processes of the form $\int_{\R} (I_+^dg)(x) dL_\alpha(x)$. The following norm is based on theorem 3.2 of Pipiras and Taqqu (2000), which in that case it comes from an inner product. Now, inspired by Proposition \ref{MRI_prop_1} and using the previous remark, we give the following definition.



\begin{defin} \label{MRI_def_1}
	For a fixed $1 < p < \alpha$ and $g \in L^1(\R) \cap L^\alpha(\R)$, let the \textit{norm} $||\cdot||_{\alpha,p}$ be defined as
	\begin{equation*} \label{MRI_eq_9}
		||g||_{\alpha,p} = \E[|L_\alpha(1)|^p]^{\frac{1}{p}} \left(\int_{\R} |(I_-^dg)(x)|^\alpha dx \right)^{\frac{1}{\alpha}} = \E[|L_\alpha(1)|^p]^{\frac{1}{p}} \|I_-^dg\|_{\alpha} \; .
	\end{equation*}
\end{defin}

\begin{prop}
	The expression given in Definition \ref{MRI_def_1} defines a norm. Moreover, for $g \in L^1(\R)\cap L^\alpha(\R)$, there exist $M,N > 0$ such that
	\begin{equation}\label{MRI_eq_10}
		||g||_{\alpha,p} \le M||g||_1 + N||g||_\alpha \; .
	\end{equation}
\end{prop}

\noindent \textbf{Proof.}
First, note that $||g||_{\alpha,p} = D||I_-^dg||_\alpha$, where $D = \left(\E[|L_\alpha(1)|^p]\right)^{\frac{1}{p}}$ does not depend on $g$ and $||\cdot||_\alpha$ is the usual $L^\alpha(\R)$ norm.
\begin{itemize}
	\item [(i)] since $D>0$, $||g||_{\alpha,p} = 0 \iff ||I_-^dg||_\alpha = 0 \iff I_-^dg = 0 \iff g = \mathcal{D}_-^d(I_-^dg) = 0$;
	\item [(ii)] for $a \in \R$, $||ag||_{\alpha,p} = D||I_-^dag||_\alpha = D||aI_-^dg||_\alpha = |a|||g||_{\alpha,p}$;
	\item [(iii)] for $f, g \in L^1(\R)\cap L^\alpha(\R), ||f+g||_{\alpha,p} = D||I_-^d(f+g)||_\alpha = D||I_-^df + I_-^dg)||_\alpha \le D||I_-^df||_\alpha + D||I_-^dg||_\alpha = ||f||_{\alpha,p} + ||g||_{\alpha,p}$ .
\end{itemize}
Finally, (\ref{MRI_eq_10}) can be seen from the proof of Proposition \ref{MRI_prop_2}: by (\ref{MRI_eq_8}) we get that $G(h) = \int_{\R} h(x)(I_-^dg)(x) dx$ is a bounded linear functional in $L^\beta(\R)$, where $\frac{1}{\alpha} + \frac{1}{\beta} = 1$. Moreover,
\begin{equation*}
	||(I_-^dg)||_\alpha = \sup_{||h||_\beta = 1} |G(h)| \le \left[\frac{1}{\Gamma(d)}\left(\frac{||g||_\alpha}{d} + \frac{||g||_1}{(\alpha(1-d) - 1)^\frac{1}{\alpha}}\right)\right] \; .
\end{equation*}
And (\ref{MRI_eq_10}) follows by setting
\begin{equation*}
	M = \frac{\left(\E[|L_\alpha(1)|^p]\right)^{\frac{1}{p}}}{\Gamma(d)}\frac{1}{(\alpha(1-d) - 1)^\frac{1}{\alpha}} > 0 \text{ and } N = \frac{\left(\E[|L_\alpha(1)|^p]\right)^{\frac{1}{p}}}{\Gamma(d)} \frac{1}{d} > 0.
\end{equation*}
\qed

\noindent We are now able to define the space which will allow us to build the desired stochastic integrals.

\begin{defin} \label{MRI_def_2}
	Let $B_{\alpha,p}$ be the space $L^1(\R) \cap L^\alpha(\R)$ under the equivalence relation $\sim_{\alpha,p}$ in which for $f,g \in L^1(\R) \cap L^\alpha(\R)$,
	\begin{equation}
		f \sim_{\alpha,p} g \iff ||f - g||_{\alpha,p} = 0 \; .
	\end{equation}
\end{defin}

\subsection{Stochastic Integral and Integration by Parts Formula}

We proceed with the construction of the stochastic integral concerning the linear fractional stable motion. Details on the stochastic integration concerning the fractional Brownian motion can be found in Pipiras and Taqqu (2000), Nualart (2006), and Carmona et al. (2001) and concerning fractional L\'{e}vy processes as a generalization of the fBm using Malliavin Calculus in He (2014).

\noindent Firstly, let us recall the LFSM defined in (\ref{I_eq_1}): let $\{L_\alpha(t)\}_{t \in \R}$ be a $S\alpha S$ process with index of stability $1 < \alpha <2$ and Lebesgue control measure $dm(x) = dx$ and let $0 < d < 1 - \frac{1}{\alpha}$ be the fractional integration parameter. Then, taking $a = (\Gamma(d+1))^{-1}$ and $b = 0$ and using the notation of Remark \ref{DPR_remark_1}, we define
\begin{equation} \label{MRI_eq_11}
	M_d(t) := \frac{1}{\Gamma(d+1)} \int_{-\infty}^{\infty} \left[(t-x)_+^d - (-x)_+^d\right] \; dL_\alpha(x) \;, \; \; t \in \R \; .
\end{equation}

\noindent Equipped with the $B_{\alpha,p}$ space given in Definition \ref{MRI_def_2} and the norm $||\cdot||_{\alpha,p}$ given in Definition \ref{MRI_def_1}, we want to define, for functions $g \in B_{\alpha,p}$, the stochastic integral
\begin{equation}
	I_{M_d}(g) = \int_{\R} g(x) dM_d(x).
\end{equation}
\noindent As usual, we start by defining the integral for simple functions. Let $\phi : \R \to \R$ be of the form
\begin{equation*}
	\phi(x) = \sum_{i=1}^{n-1}a_i \I_{(t_i,t_{i+1}]}(x),
\end{equation*}
where $a_i \in \R$ and $-\infty < t_1 < t_2 < \cdots < t_n < \infty$. It is easy to see that $\phi \in B_{\alpha,p}$. So we define the integral
\begin{equation*}
	I_{M_d}(\phi) = \int_{\R} \phi(x) dM_d(x) = \sum_{i=1}^{n-1} a_i[M_d(t_{i+1}) - M_d(t_i)].
\end{equation*}
\noindent Note that the integral $I_{M_d}(\cdot)$ is linear for simple functions. Moreover, we have the proposition below.

\begin{prop} \label{MRI_prop_3}
	Let $\phi : \R \to \R$ be a simple function. Then
	\begin{equation} \label{MRI_eq_12}
		\int_{\R} \phi(x) dM_d(x) = \int_{\R} (I_-^d\phi)(x) \; dL_\alpha(x)
	\end{equation}
	and $\phi \to I_{M_d}(\phi)$ is an isometry between $B_{\alpha,p}$ and $L^p(\Omega, \mathbb{P})$, in the sense that $||I_{M_d}(\phi)||_{L^p(\Omega, \mathbb{P})}^p$ $ = ||\phi||_{\alpha,p}^p$, where $(\Omega, \mathcal{F}, \mathbb{P})$ is the underlying probability space.
\end{prop}

\noindent \textbf{Proof.}
	It suffices to show (\ref{MRI_eq_12}) for indicator functions $\phi(x) = \I_{(0,t]}(x)$, $t > 0$, because every other simple function is a linear combination of such functions. In fact,
	\begin{equation*}
		\int_{\R} \phi(x) dM_d(x) = \int_{\R} \I_{(0,t]}(x) dM_d(x) = M_d(t)
	\end{equation*}
	\noindent and for the right-hand side of (\ref{MRI_eq_12}), we have
	\begin{align}
			\int_{\R} (I_-^d\phi)(x) \; dL_\alpha(x) &= \frac{1}{\Gamma(d)}\int_{\R}\int_{x}^{\infty} (s-x)^{d-1}\I_{(0,t]}(s) dsdL_\alpha(x) \nonumber\\
			&= \frac{1}{\Gamma(d+1)} \int_{\R}[(t-x)_+^d - (-x)_+^d] dL_\alpha(x) = M_d(t).
	\end{align}
	\noindent Furthermore, for all simple functions $\phi$ it follows from Proposition \ref{MRI_prop_1} that
\begin{align} \label{MRI_eq_13}
	\nonumber
	||I_{M_d}(\phi)||_{L^p(\Omega, \mathbb{P})}^p &= \E\left[\left| \int_{\R} (I_-^d\phi)(x) dL_\alpha(x) \right|^p\right] = \E[|L_\alpha(1)|^p] \left(\int_{\R} |(I_-^d\phi)(x)|^\alpha dx\right)^\frac{p}{\alpha}  \\
	&= ||\phi||_{\alpha,p}^p.
\end{align}
\vspace{-0.6cm}\qed
\noindent Having the integral $I_{M_d}(\cdot)$ well defined for simple functions, we are ready to prove the next theorem and to define the integral for any $g \in B_{\alpha,p}$.

\begin{teo} \label{MRI_teo_1}
	Let $M_d = \{M_d(t)\}_{t \in \R}$ be the linear fractional stable motion defined in (\ref{MRI_eq_11}) and let the function $g : \R \to \R$ be in $B_{\alpha,p}$. Then there exists a sequence of simple functions $\phi_n : \R \to \R$, $n \in \N$, satisfying $||\phi_n - g||_{\alpha,p} \to 0$ as $n \to \infty$ such that $I_{M_d}(\phi_n)$ converges in $L^p(\Omega, \mathbb{P})$ to a limit denoted by $I_{M_d}(g) = \int_{\R} g(x)dM_d(x)$ and $I_{M_d}(g)$ is independent of the approximating sequence $\phi_n$. Moreover, the isometry property holds
	\begin{equation} \label{MRI_eq_14}
	||I_{M_d}(g)||_{L^p(\Omega, \mathbb{P})} = ||g||_{\alpha,p} \; .
	\end{equation}
\end{teo}

\noindent \textbf{Proof.}
	As we know, simple functions are dense in $L^1(\R)\cap L^\alpha(\R)$ and $L^1(\R)\cap L^\alpha(\R)$ is dense in $B_{\alpha,p}$ by construction of the space $B_{\alpha,p}$. Also (\ref{MRI_eq_10}) holds, so whenever a sequence converges in $L^1(\R)\cap L^\alpha(\R)$ it will converge in $B_{\alpha,p}$ as well. Therefore, the simple functions are dense in $B_{\alpha,p}$. Hence, there exists a sequence $\phi_n$ of simple functions such that $||\phi_n - g||_{\alpha,p} \to 0$, as $n \to \infty$. It follows from the isometry property (\ref{MRI_eq_13}) that $\int_{\R} \phi_n(x) dM_d(x)$ converges in $L^p(\Omega, \mathbb{P})$ towards a limit denoted by $\int_{\R} g(x) dM_d(x)$ and the isometry property is preserved in this procedure, so (\ref{MRI_eq_14}) holds. Finally, (\ref{MRI_eq_14}) implies that the limit denoted by the integral $\int_{\R} g(x) dM_d(x)$ is the same for all approximating sequences $\phi_n$ converging to $g$. \vspace{-0.4cm}\qed

\noindent Note that convergence in the $L^p$ sense implies convergence in probability, which is enough to define the integrals. Usually, they are defined with equality in probability or even in the $L^p$ sense, most commonly in $L^2$ (see Kuo, 2006). The next theorem states that we can interchange the fractional integration between the integrand function and the process concerning which we are integrating. This formula reminds the well-known integration by parts formula from traditional calculus.

\begin{teo} \label{MRI_teo_2}
	Let $g \in B_{\alpha,p}$. Then
	\begin{equation}
		\int_{\R} g(x) dM_d(x) = \int_{\R} (I_-^dg)(x) \; dL_\alpha(x) \; .
	\end{equation}
\end{teo}

\noindent \textbf{Proof.}
	It follows from Proposition \ref{MRI_prop_3} and Theorem \ref{MRI_teo_1}. \vspace{-0.4cm}\qed

\noindent Now, as a direct application of the previous results, we state a convergence theorem that focuses on using the isometry property (\ref{MRI_eq_14}) and the norm $||\cdot||_{\alpha,p}$ to make sense the analysis of the convergence in $L^p(\Omega, \mathbb{P})$ of $I_{M_d}(g_n) \to I_{M_d}(g)$, as $n \to \infty$, for a suitable choice of $g_n$ and $g$. We will study convergence on $L^p(\Omega, \mathbb{P})$ by dealing with convergence in $L^1(\R)\cap L^\alpha(\R)$.

\noindent Consider $g,g_n$ real measurable functions in $B_{\alpha,p}$ such that the stochastic integrals $I_{M_d}(g_n)$ and $I_{M_d}(g)$ are all well defined according to Theorem \ref{MRI_teo_1}. We say that the sequence of processes $I_{M_d}(g_n)$ converges to the limit process $I_{M_d}(g)$ if
\vspace{.3cm}
\begin{equation}
	||I_{M_d}(g_n) - I_{M_d}(g)||_{L^p(\Omega, \mathbb{P})} \to 0
\end{equation}
\noindent and we denote
\begin{equation}
	\lim\limits_{n \to \infty} \int_{\R} g_n(x) dM_d(x) = \int_{\R} g(x) dM_d(x) \; ,
\end{equation}
\noindent with equality in the $L^p(\Omega, \mathbb{P})$ sense. Of course, it also implies convergence in probability for $\mathbb{P}$. Equality has the same meaning as in Theorem \ref{MRI_teo_1}.

\begin{teo} \label{MRI_teo_3}
	If $g_n \to g$ in both $L^1(\R)$ and $L^\alpha(\R)$, then $I_{M_d}(g_n) \to I_{M_d}(g)$ in $L^p(\Omega, \mathbb{P})$ and
	\begin{equation}
		\lim\limits_{n \to \infty} \int_{\R} g_n(x) dM_d(x) = \int_{\R} g(x) dM_d(x).
	\end{equation}
\end{teo}

\noindent \textbf{Proof.}
	By (\ref{MRI_eq_14}) and (\ref{MRI_eq_10}) we have
	\begin{equation}
		||I_{M_d}(g_n) - I_{M_d}(g)||_{L^p(\Omega, \mathbb{P})} = ||g_n - g||_{\alpha,p} \le M||g_n - g||_1 + N||g_n - g||_\alpha
	\end{equation}
	and it follows that $||I_{M_d}(g_n) - I_{M_d}(g)||_{L^p(\Omega, \mathbb{P})} \to 0$, as $n \to \infty$.  \qed

\section{Main Results II}
\renewcommand{\theequation}{\thesection.\arabic{equation}}
\setcounter{equation}{0}

In this section, we are interested in defining a class of continuous-time moving average stable processes that present the property of long-range dependence. Let $\{L_\alpha(t)\}_{t \in \R}$ be a $S\alpha S$ process with index of stability  $1 < \alpha < 2$ and Lebesgue control measure. Accordingly to Samorodnitsky and Taqqu (1994), if $g$ is a real measurable function defined on $\R$ satisfying $\int_{\R} |g(x)|^\alpha dx < \infty$, and define
\begin{equation}
X(t) = \int_{-\infty}^{\infty} g(t - x) dL_\alpha(x), \; \; t \in \R \; ,
\end{equation}
\noindent then $\{X(t)\}_{t \in \R}$ is well defined and stationary. A process of this form is called \itshape $S\alpha S$ moving average process. \normalfont  The function $g$ is often called the \emph{kernel function}. 

\noindent It's known in the literature that one can obtain a long-memory process by considering moving average processes whose kernel function is the fractional integral of a \textit{short-memory} function (see Samorodnitsky, 2006 and Marquardt, 2006). Henceforth, we are considering functions with \textit{short memory} and with positive support, which means that we will assume the following conditions:

\begin{itemize}
	\item[(C1)] $g(t) = 0$, for all $t < 0$;
	\item[(C2)] $|g(t)| \leq C e^{-ct}$, for some constants $C>0$ and $c>0$.
\end{itemize}
\noindent Note that every such function $g$ satisfying (C1) and (C2) belongs to $B_{\alpha,p}$.

\noindent Let us recall the left-sided Riemann-Liouville fractional integral of order $d \in (0,1)$ of a function $g \in L^1(\R)$, as defined in \eqref{rlRL}:
\begin{equation*}
	(I_+^dg)(x) = \frac{1}{\Gamma(d)} \int_{-\infty}^{x} g(t)(x-t)^{d-1} dt = \frac{1}{\Gamma(d)} \int_{0}^{\infty} g(x - s)(s)^{d-1} ds.
\end{equation*}
\noindent So for $1 < \alpha < 2$ and $1 < p < \alpha$ the idea is to define the moving average stable process using the left-sided Riemann-Liouville fractional integral of a kernel function $g$ satisfying conditions (C1) and (C2).

\begin{defin} [FIMA Stable process]
	Let $\{L_\alpha(t)\}_{t \in \R}$ be a $S\alpha S$ process with index of stability  $1 < \alpha < 2$ and Lebesgue control measure. For $d \in (0,1 - \frac{1}{\alpha})$, define
	\begin{equation} \label{MR2_eq_1}
		Y_d(t) = \int_{-\infty}^{t} (I_+^dg)(t - x) dL_\alpha(x) \; , t \in \R \; .
	\end{equation}
\end{defin}

\begin{teo}
	The FIMA Stable process $\{Y_d(t)\}_{t \in \R}$ in (\ref{MR2_eq_1}) is well defined and stationary.
\end{teo}
\noindent \textbf{Proof.}
	We know that $g \in B_{\alpha,p}$ implies $g \in L^1(\R)\cap L^\alpha(\R)$, so by Remark \ref{MRI_remark_1} it follows that $\int_{\R} |(I_+^dg)(x)|^\alpha dx < \infty$ and the process is well defined. Now, for $t_1, \cdots, t_d, h, \theta_1, \cdots, \theta_d \in \mathbb{R}$, setting $y = x - h$, we have
	\begin{align}
		\left|\left|\sum_{i=1}^{d} \theta_iY_d(t_i + h) \right|\right|_\alpha^\alpha
		&= \int_{-\infty}^{\infty} \left|\sum_{i=1}^{d} \theta_i \I_{(-\infty, t_i + h]}(x) (I_+^dg)(t_i + h - x)\right|^\alpha dx \nonumber\\
		&= \int_{-\infty}^{\infty} \left|\sum_{i=1}^{d} \theta_i \I_{(-\infty, t_i]}(y) (I_+^dg)(t_i - y)\right|^\alpha dy = \left|\left|\sum_{i=1}^{d} \theta_iY_d(t_i) \right|\right|_\alpha^\alpha.
	\end{align} \qed

\noindent Now, using Theorem \ref{MRI_teo_2}, we can represent the FIMA Stable process with the linear fractional stable motion.

\begin{teo} \label{MRII_teo_2}
	Let $\{Y_d(t)\}_{t \in \R}$ be the FIMA Stable process defined in (\ref{MR2_eq_1}) and $\{M_d(t)\}_{t \in \R}$ be the linear fractional stable motion defined in (\ref{MRI_eq_11}). Then, $Y_d$ can be represented as
	\begin{equation} \label{MR2_eq_2}
		Y_d(t) = \int_{-\infty}^{t} g(t - x) dM_d(x) \; , \; \; t \in \R.
	\end{equation}
\end{teo}

\noindent \textbf{Proof.} For each $t \in \R$, consider $g_t(x) := g(t-x)$, so from Theorem \ref{MRI_teo_2}, we have
	\begin{equation*}
		\int_{-\infty}^{t} g(t - x) dM_d(x) = \int_{-\infty}^{t} g_t(x) dM_d(x) = \int_{-\infty}^{t} (I_-^dg_t)(x) dL_\alpha(x).
	\end{equation*}
	By using the definition of the right-sided fractional integral $(I_-^dg_t)(x)$ and setting $w = u - x$ we get
	\begin{align} \label{MRII_eq_3}
\nonumber
			\int_{-\infty}^{t} (I_-^dg_t)(x) dL_\alpha(x)
			&= \frac{1}{\Gamma(d)} \int_{-\infty}^{t} \int_{x}^{\infty} (u - x)^{d-1} g_t(u) du \; dL_\alpha(x) \\
\nonumber
			&= \frac{1}{\Gamma(d)} \int_{-\infty}^{t} \int_{x}^{\infty} (u - x)^{d-1} g(t - u) du \; dL_\alpha(x) \\
			&= \frac{1}{\Gamma(d)} \int_{-\infty}^{t} \int_{0}^{\infty} w^{d-1} g(t - x - w) dw \; dL_\alpha(x).
	\end{align}
	
	\noindent Finally, using the definition of the left-sided fractional integral $(I_+^dg)(t - x)$ we get the result
	\begin{equation} \label{MRII_eq_4}
		\frac{1}{\Gamma(d)} \int_{-\infty}^{t} \int_{0}^{\infty} w^{d-1} g(t - x - w) dw \; dL_\alpha(x) =
		\int_{-\infty}^{t} (I_+^dg)(t - x) dL_\alpha(x) = Y_d(t) \; .
	\end{equation} \qed

\begin{obs}
	Note that this representation allows us to work with a rapid decay kernel function $g$ instead of a slow decay kernel function $(I_+^dg)$, which can be useful for more efficient simulation algorithms. It also can be used to show a convergence of sequences of such processes by dealing with the associated moving average kernel function. See Theorems \ref{MRI_teo_3} and \ref{TheoMRII4.5}.
\end{obs}

\noindent Now we turn our attention to the long-range dependence, or long-memory property of the process $\{Y_d(t)\}_{t \in \R}$. As we know, for stationary Gaussian processes one can describe the long-range dependence property as the slow decaying of its autocovariance function. Since we cannot use such a function in our setting, as Maejima and Yamamoto (2003) we are going to use the following and proceed with similar calculations as the cited authors. For $\theta_1,\theta_2 \in \R$, let
\begin{align} \label{MRII_eq_1}
\nonumber
	r(t)
	&:= r(\theta_1,\theta_2; t) = \E\Big[\exp\{i(\theta_1Y_d(t) + \theta_2Y_d(0))\}\Big] \\
    &- \E\Big[\exp\{i\theta_1Y_d(t)\}\Big]\E\Big[\exp\{i\theta_2Y_d(0)\}\Big]	
\end{align}
\noindent and
\begin{align}\label{I}
\nonumber	
	I(t)
	&:= I(\theta_1,\theta_2; t) = -\log \bigg(\E\Big[\exp\{i(\theta_1Y_d(t) + \theta_2Y_d(0))\}\Big]\bigg) \nonumber\\
	&+ \log \bigg(\E\Big[\exp\{i\theta_1Y_d(t)\}\Big]\bigg) + \log \bigg(\E\Big[\exp\{i\theta_2Y_d(0)\}\Big]\bigg) \; .
\end{align}

\noindent Setting $$K(\theta_1,\theta_2;t) = \E\big[\exp\{i\theta_1Y_d(t)\}\big]\E\big[\exp\{i\theta_2Y_d(0)\}\big]$$ $$= \E\big[\exp\{i\theta_1Y_d(0)\}\big]\E\big[\exp\{i\theta_2Y_d(0)\}\big] =: K$$ which does not depend on $t$ since $\{Y_d(t)\}_{t \in \R}$ is stationary. There is a relationship among $r$, $I$ and $K$ given by

\begin{equation}\label{r_I_K}
	r(\theta_1,\theta_2;t) = K(\theta_1,\theta_2;t)\big(e^{-I(\theta_1,\theta_2;t)} - 1\big).
\end{equation}

\noindent Furthermore, if $I(t) \to 0$ as $t \to \infty$, then $r(t) \sim -KI(t)$ as $t \to \infty$, which means that the quantities are asymptotically equal. Note that the quantity $r(\theta_1,\theta_2;t)$ is the difference between the joint characteristic function of the pair $(Y_d(t), Y_d(0))$ and the product of the marginal characteristic functions of $Y_d(t)$ and $Y_d(0)$. So if the process is independent, $r(\theta_1,\theta_2;t) = 0$, and if the process is Gaussian or has finite second-order moments, the quantity $-I(1,-1;t)$ coincides with the autocovariance function. The quantity $r(t)$ has been proved to be a proper tool for describing the dependence structure of stable processes and has been used by several authors such as Astrauskas et al. (1991); L\'evy and Taqqu (1991); Magdziarz and Weron (2007); Magdziarz (2007). Besides, it is used in Maruyama (1970) relating it to the mixing property for stationary infinitely divisible processes. Also, Samorodnitsky and Taqqu (1994) point out that the quantity $I(\cdot)$ can be used to show that two processes are different by showing that their respective quantities are asymptotically different. Let us then define what we mean by the long-range dependence or long-memory property in our setting.

\begin{defin} \label{MRII_def_1}
	Let $\{X(t)\}_{t \in \R}$ be a stationary $S\alpha S$ process and $r(t)$ as in (\ref{MRII_eq_1}). We say the process $\{X(t)\}_{t \in \R}$ has \textit{long-range dependence} or \textit{long-memory property} if there exists $\theta \in (0,1)$ and a constant $c_r > 0$ such that $r(t)$ satisfies
	\begin{equation}
		\lim\limits_{t \to \infty} |r(t)|\frac{t^\theta}{c_r} = 1 \; ,
	\end{equation}
	or, equivalently,
	\begin{equation*}
		|r(t)| \sim c_r t^{-\theta} \; .
	\end{equation*}
\end{defin}

\noindent Note that since the process $\{Y_d(t)\}_{t \in \R}$ is stationary,
\begin{equation*}
	(Y_d(t), Y_d(0)) \overset{d}{=} (Y_d(0), Y_d(-t)) \; ,
\end{equation*}

\noindent thus $r(\theta_1,\theta_2;t) = r(\theta_2,\theta_1;-t)$. Therefore, the asymptotic behavior of $r(\theta_1,\theta_2;t)$ as $t \to \infty$ is essentially the same as $r(\theta_1,\theta_2;t)$ as $t \to -\infty$. Since the computations for $t<0$ as $t \to -\infty$ are a bit easier in our setting, as in Maejima and Yamamoto (2003), we will analyze the asymptotic behavior of $r(\theta_1,\theta_2; t)$ for $t<0$, as $t \to -\infty$.

\noindent Before we give the proof of the asymptotic behavior of $r(\theta_1, \theta_2; t)$, for $t<0$, as $t\to -\infty$, for a FIMA process (see Theorem 4.3 below), we shall give three auxiliary lemmas.

\begin{lema} \label{MRII_lema_1}
	For $1 < \alpha < 2$ and for $r,s \in \R$, it is true that
	\begin{equation}
		||r + s|^\alpha - |r|^\alpha - |s|^\alpha| \le \alpha|r||s|^{\alpha-1} + (\alpha+1)|r|^\alpha.
	\end{equation}
\end{lema}

\noindent \textbf{Proof.} See lemma 5.2 in Maejima and Yamamoto (2003).\vspace{-0.2cm}\qed

We define the following functions $\phi(\cdot,x)$ and $\psi(\cdot,x)$,
for $x >1$,
\begin{equation}\label{phipsi}
	\phi(t,x) = \frac{1}{\Gamma(d)} \int_{1}^{x} \theta_1 g(t(1-u))(x-u)^{d-1} du \ \ \mbox{and} \ \
    \psi(t,x) = \frac{1}{\Gamma(d)} \int_{0}^{x} \theta_2 g(-tu)(x-u)^{d-1} du.
\end{equation}

\begin{lema} \label{MRII_lema_2}
	Let the functions $\phi(\cdot,x)$ and $\psi(\cdot,x)$ be defined in \eqref{phipsi}. Then, there exist constants $K_1,K_2>0$ such that, for any $t < 0$ and $x > 1$ we have
	\begin{equation} \label{MRII_eq_7}
		|t\phi(t,x)| \le K_1(x-1)^{d-1} \ \ \mbox{and } \ \
	|t\psi(t,x)| \le K_2x^{d-1}.
	\end{equation}
\end{lema}

\noindent \textbf{Proof.} First, note that setting $w = u+1$, we have
\begin{align}
\nonumber
\theta_2^{-1} \psi(t,x-1) &= \frac{1}{\Gamma(d)} \int_{0}^{x-1} g(-tu)(x-1-u)^{d-1} du \\
& =\frac{1}{\Gamma(d)} \int_{1}^{x} g(t(1-w))(x-w)^{d-1} dw  = \theta_1^{-1} \phi(t,x).
\end{align}
\noindent Thus it suffices to show (\ref{MRII_eq_7}) for $\phi(\cdot,x)$. In fact, we have
	\begin{align}
		&|t\phi(t,x)| = \left|\frac{1}{\Gamma(d)} \int_{1}^{x} t\theta_1 g(t(1-u))(x-u)^{d-1} du\right| \nonumber\\
		&\le \left|\frac{1}{\Gamma(d)} \int_{1}^{\frac{x+1}{2}} t\theta_1 g(t(1-u))(x-u)^{d-1} du\right| + \left|\frac{1}{\Gamma(d)} \int_{\frac{x+1}{2}}^{x} t\theta_1 g(t(1-u))(x-u)^{d-1} du\right| \nonumber\\
		& = I_1 + I_2.
	\end{align}
	Recall that $t<0$. For $I_1$, it follows that
	\begin{equation} \label{MRII_eq_9}
		I_1 \le \frac{|t\theta_1 C|}{\Gamma(d)} \left(\frac{x-1}{2}\right)^{d-1} \int_{1}^{\frac{x+1}{2}} e^{ct(u-1)} du \le \frac{|\theta_1C|}{\Gamma(d)c} \left(\frac{x-1}{2}\right)^{d-1}
	\end{equation}
	\noindent and
	\begin{align} \label{MRII_eq_10}
\nonumber
		&I_2 \le \frac{|t\theta_1 C|}{\Gamma(d)} e^{\frac{ct(x-1)}{2}} \int_{\frac{x+1}{2}}^{x} (x-u)^{d-1} du = \frac{|\theta_1 C|}{\Gamma(d)} \left(\frac{x-1}{2}\right)^d |t|e^{\frac{ct(x-1)}{2}} \\
& \le \frac{|\theta_1 C|}{\Gamma(d)ce} \left(\frac{x-1}{2}\right)^{d-1},
	\end{align}
	since $z(t) := |t|e^{\frac{ct(x-1)}{2}}$ takes its maximum value $\frac{2}{c(x-1)e}$ at $t = \frac{-2}{c(x-1)}$.
	\noindent Combining (\ref{MRII_eq_9}) and (\ref{MRII_eq_10}), we have
	\begin{equation}
		|t\phi(t,x)| \le I_1 + I_2 \le \frac{|\theta_1C|}{\Gamma(d)c} \left(\frac{x-1}{2}\right)^{d-1} + \frac{|\theta_1 C|}{\Gamma(d)ce} \left(\frac{x-1}{2}\right)^{d-1} \le K_1(x-1)^{d-1} \; ,
	\end{equation}
\noindent which is (\ref{MRII_eq_7}). \vspace{-0.3cm}\qed

\begin{lema} \label{MRII_lema_3}
Let the functions $\phi(\cdot,x)$ and $\psi(\cdot,x)$ be defined in \eqref{phipsi}.
	There exists a constant $L \in \R$ such that
	\begin{equation} \label{MRII_eq_8}
		\lim\limits_{t \to -\infty} t\phi(t,x) = L\theta_1 (x-1)^{d-1} \ \ \mbox{and} \ \
	     \lim\limits_{t \to -\infty} t\psi(t,x) = L\theta_2 x^{d-1} \; .
	\end{equation}
\end{lema}

\noindent \textbf{Proof.} For (\ref{MRII_eq_8}) it suffices to show for the function $\phi(\cdot,x)$. It follows that,
	\begin{align} \label{MRII_eq_11}
		\lim\limits_{t \to -\infty} t\phi(t,x)
		&= \lim\limits_{t \to -\infty} \frac{1}{\Gamma(d)} \int_{1}^{\frac{x+1}{2}} t\theta_1 g(t(1-u))(x-u)^{d-1} du \nonumber\\
		&+ \lim\limits_{t \to -\infty} \frac{1}{\Gamma(d)} \int_{\frac{x+1}{2}}^{x} t\theta_1 g(t(1-u))(x-u)^{d-1} du = J_1 + J_2.
	\end{align}
	So, set $w = t(u-1)$ to get
	\begin{align}
			J_1
			&= \lim\limits_{t \to -\infty} \frac{1}{\Gamma(d)} \int_{0}^{\frac{t(x-1)}{2}} \theta_1 g(-w)\left(x-\frac{w}{t} -1\right)^{d-1} dw \nonumber\\
			&= -\lim\limits_{t \to -\infty} \frac{\theta_1}{\Gamma(d)} \int_{-\infty}^{0} \I_{\left[\frac{t(x-1)}{2}, 0\right]}(w) g(-w)\left(x-\frac{w}{t} -1\right)^{d-1} dw
	\end{align}
	but by condition (C2),
	\begin{align}
\nonumber
		&\left|\I_{\left[\frac{t(x-1)}{2}, 0\right]}(v) g(-w)\left(x-\frac{w}{t} -1\right)^{d-1} \right| \le C\left|\I_{\left[\frac{t(x-1)}{2}, 0\right]}(w) e^{cv}\left(x-\frac{w}{t} -1\right)^{d-1} \right| \\
        &\le Ce^{cw}\left(\frac{x-1}{2}\right)^{d-1},
	\end{align}
	\noindent which is $w$ integrable in $(-\infty,0)$. Thus, by the Dominated Convergence theorem and setting $u = -w$, we get
	\begin{equation}
		J_1 = -\frac{\theta_1}{\Gamma(d)} \int_{-\infty}^{0} g(-w)\left(x-1\right)^{d-1} dw = \left[\frac{1}{\Gamma(d)} \int_{0}^{\infty} g(u)du\right] \theta_1 (x-1)^{d-1} \; .
	\end{equation}
	And for $J_2$,
	\begin{align} \label{MRII_eq_12}
\nonumber
		|J_2| &\le \lim\limits_{t \to -\infty} \frac{|\theta_1|}{\Gamma(d)} \int_{\frac{x+1}{2}}^{x} |tg(t(1-u))|(x-u)^{d-1} du \\
		&\le \lim\limits_{t \to -\infty} \frac{|t\theta_1|}{\Gamma(d)} \int_{\frac{x+1}{2}}^{x} e^{ct(u-1)}(x-u)^{d-1} du \le \lim\limits_{t \to -\infty} \frac{|t\theta_1|}{\Gamma(d)} e^{\frac{ct(x-1)}{2}} \int_{\frac{x+1}{2}}^{x}(x-u)^{d-1} du = 0.
	\end{align}
	\noindent Therefore, from equations (\ref{MRII_eq_10})-(\ref{MRII_eq_11}), it follows that
	\begin{equation}
		\lim\limits_{t \to -\infty} t\phi(t,x) = \left[\frac{1}{\Gamma(d)} \int_{0}^{\infty} g(u)du\right] \theta_1 (x-1)^{d-1}
	\end{equation}
	\noindent and we get (\ref{MRII_eq_8}) by setting
	\begin{equation*}
		L = \frac{1}{\Gamma(d)} \int_{0}^{\infty} g(u)du \; .
	\end{equation*}
\vspace{-0.3cm}\qed


\noindent Now we are ready for Theorem \ref{MRII_teo_1} below. Recall that for a $S\alpha S$ process $\{X(t)\}_{t \in \R}$ given by (\ref{MRI_eq_1}), the Proposition \ref{DPR_prop_1} states that
\begin{equation} \label{MRII_eq_2}
	\E\left[\exp\left\{ i\theta \int_{\R} f_t(x) dL_\alpha(x) \right\}\right] = \exp\left\{ -|\theta|^\alpha \int_{\R} |f_t(x)|^\alpha dx \right\} \; .
\end{equation}

\begin{teo} \label{MRII_teo_1}
	Let $\{Y_d(t)\}_{t \in \R}$ be the FIMA Stable process defined in (\ref{MR2_eq_1}). If $1 < \alpha < 2$ and $0 < d < 1 - \frac{1}{\alpha}$, then, as $t \to -\infty$,
	\begin{equation}
		r(t) = r(\theta_1,\theta_2;t) \sim -KC|t|^{\alpha(d-1)+1} \; ,
	\end{equation}
	where
	\begin{equation}
		K = \exp\Bigg\{-\frac{\theta_1^\alpha + \theta_2^\alpha}{\Gamma(d)^\alpha} \int_{-\infty}^{0}\left| \int_{x}^{0} g(-u)(u-x)^{d-1} du\right|^\alpha dx \Bigg\} \; ,
	\end{equation}
	and
	\begin{equation} \label{MRII_eq_6}
		C\hspace{-0.1cm} = \left(\frac{1}{\Gamma(d)} \int_{0}^{\infty} g(u)du\right)^\alpha
		\int_{1}^{\infty} \Big(\left| \theta_1(x-1)^{d-1} + \theta_2x^{d-1} \right|^\alpha - \left| \theta_1(x-1)^{d-1}\right|^\alpha - \left|\theta_2x^{d-1} \right|^\alpha\Big) dx.
	\end{equation}
\end{teo}

\noindent \textbf{Proof.} We start by calculating $K=\E\big[\exp\{i\theta_1Y_d(0)\}\big]\E\big[\exp\{i\theta_2Y_d(0)\}\big]$. By \eqref{MRII_eq_2}, (\ref{MRII_eq_3}) and (\ref{MRII_eq_4})
\begin{align}
\nonumber
	K \hspace{-0.1cm}&=\hspace{-0.1cm}\exp\Bigg\{-\frac{1}{\Gamma(d)^\alpha} \int_{-\infty}^{0} \bigg(\left| \int_{x}^{0} \theta_1 g(-u)(u-x)^{d-1} du\right|^\alpha - \left|\int_{x}^{0} \theta_2 g(-u)(u-x)^{d-1} du\right|^\alpha \bigg) dx \Bigg\} \nonumber\\
	&= \exp\Bigg\{-\frac{\theta_1^\alpha + \theta_2^\alpha}{\Gamma(d)^\alpha} \int_{-\infty}^{0}\left| \int_{x}^{0} g(-u)(u-x)^{d-1} du\right|^\alpha dx \Bigg\} \; .
\end{align}

\noindent Now, let us compute $I(\cdot)$. Again by (\ref{MRII_eq_2}), (\ref{MRII_eq_3}) and (\ref{MRII_eq_4}), we have
\begin{align}
\nonumber
	I(t) &= -\log \bigg(\E\Big[\exp\{i(\theta_1Y_d(t) + \theta_2Y_d(0))\}\Big]\bigg) + \log \bigg(\E\Big[\exp\{i\theta_1Y_d(t)\}\Big]\bigg) \\
\nonumber
         & + \log \bigg(\E\Big[\exp\{i\theta_2Y_d(0)\}\Big]\bigg)\\
\nonumber
	&=\frac{1}{\Gamma(d)^\alpha} \int_{-\infty}^{t} \bigg(\left| \int_{x}^{t} \theta_1 g(t-u)(u-x)^{d-1} du + \int_{x}^{0} \theta_2 g(-u)(u-x)^{d-1} du\right|^\alpha \nonumber\\
	&- \left|\int_{x}^{t} \theta_1 g(t-u)(u-x)^{d-1} du\right|^\alpha - \left|\int_{x}^{0} \theta_2 g(-u)(u-x)^{d-1} du\right|^\alpha \bigg) dx.
\end{align}

\noindent Set $x = xt$ and $u = ut$ and recall that $t<0$ to get
\begin{align}
	I(t)
	&=\frac{1}{\Gamma(d)^\alpha} \int_{1}^{\infty} \bigg(\bigg| \int_{1}^{x} \theta_1 g(t(1-u))(x-u)^{d-1}|t|^{d-1} |t|du \nonumber\\
	&+ \int_{0}^{x} \theta_2 g(-tu)(x-u)^{d-1}|t|^{d-1} |t|du\bigg|^\alpha - \left|\int_{1}^{x} \theta_1 g(t(1-u))(x-u)^{d-1}|t|^{d-1} |t|du\right|^\alpha \nonumber\\
	&- \left|\int_{0}^{x} \theta_2 g(-tu)(x-u)^{d-1}|t|^{d-1} |t|du\right|^\alpha \bigg) |t|dx \; .
\end{align}
\noindent Since $\left(|t|^{d-1} |t|\right)^\alpha |t| = |t|^{\alpha d + 1}$, it follows that

\begin{align}
		I(t)
		&=\frac{|t|^{\alpha d + 1}}{\Gamma(d)^\alpha} \int_{1}^{\infty} \bigg(\bigg| \int_{1}^{x} \theta_1 g(t(1-u))(x-u)^{d-1} du + \int_{0}^{x} \theta_2 g(-tu)(x-u)^{d-1} du\bigg|^\alpha \nonumber\\
		&- \left|\int_{1}^{x} \theta_1 g(t(1-u))(x-u)^{d-1} du\right|^\alpha - \left|\int_{0}^{x} \theta_2 g(-tu)(x-u)^{d-1} du\right|^\alpha \bigg) dx \nonumber\\
		&=: |t|^{\alpha d + 1} \int_{1}^{\infty} \left| \phi(t,x) + \psi(t,x) \right|^\alpha - \left| \phi(t,x)\right|^\alpha - \left|\psi(t,x) \right|^\alpha dx \; ,
\end{align}
where $\phi(\cdot, x)$ and $\psi(\cdot, x)$, for $x > 1$, are given in \eqref{phipsi}.

Therefore,
\begin{align} \label{MRII_eq_5}
	& \lim\limits_{t \to -\infty} |t|^{\alpha-\alpha d-1} I(t) =
		\lim\limits_{t \to -\infty} \int_{1}^{\infty} \left(\left| t\phi(t,x) + t\psi(t,x) \right|^\alpha - \left| t\phi(t,x)\right|^\alpha - \left|t\psi(t,x) \right|^\alpha\right) dx.	
\end{align}

\noindent We want to use the Dominated Convergence theorem to deal with the above limit.
From Lemmas \ref{MRII_lema_1}, with $r = t\psi(t,x)$ and $s = t\phi(t,x)$, and Lemma \ref{MRII_lema_2}, it follows that

\begin{align}
\nonumber
	  & \left|s+r \right|^\alpha - \left| s\right|^\alpha - \left|r \right|^\alpha \le \alpha|r||s|^{\alpha-1} + (\alpha + 1) |r|^\alpha \\
\nonumber
		&\le \alpha K_2 x^{d-1}(K_1(x-1)^{d-1})^{\alpha-1} + (\alpha+1)(K_2 x^{d-1})^\alpha \\
		&= \alpha K_1^{(\alpha-1)(d-1)}K_2x^{d-1}(x-1)^{(\alpha-1)(d-1)} + (\alpha+1)K_2^\alpha x^{\alpha(d-1)}.
\end{align}

\noindent Note that $1<\alpha<2$ implies $0<\alpha-1<1$ and $0<d<1 - \frac{1}{\alpha}$ implies $-1<d-1<-\frac{1}{\alpha}$, so $-1 < (\alpha-1)(d-1)$ and $\alpha(d-1) < -1$. Thus, $x^{\alpha(d-1)}$ clearly is in $L^1(1,\infty)$. As for $x^{d-1}(x-1)^{(\alpha-1)(d-1)}$, the integral may diverge when $x \to 1$ and $x \to \infty$. But when $x \to 1$, $x^{d-1}(x-1)^{(\alpha-1)(d-1)}$ behaves like $(x-1)^{(\alpha-1)(d-1)}$, which is integrable since $-1 < (\alpha-1)(d-1)$. And when $x \to \infty$, $x^{d-1}(x-1)^{(\alpha-1)(d-1)}$ behaves like $x^{\alpha(d-1)}$. Hence,
\begin{equation}
	\int_{1}^{\infty} \left|\alpha K_1^{(\alpha-1)(d-1)}K_2x^{d-1}(x-1)^{(\alpha-1)(d-1)} + (\alpha+1)K_2^\alpha x^{\alpha(d-1)}\right| dx < \infty.
\end{equation}

\noindent Now by applying the Dominated Convergence theorem in (\ref{MRII_eq_5}) altogether with Lemma \ref{MRII_lema_3} we obtain

\begin{align}\label{eq4_33}
\nonumber
	\noindent & \lim\limits_{t \to -\infty} |t|^{\alpha-\alpha d-1} I(t) =
	 \int_{1}^{\infty} \lim\limits_{t \to -\infty} \left| t\phi(t,x) + t\psi(t,x) \right|^\alpha - \left| t\phi(t,x)\right|^\alpha - \left|t\psi(t,x) \right|^\alpha dx \\
	&=L^\alpha \int_{1}^{\infty} \left| \theta_1(x-1)^{d-1} + \theta_2x^{d-1} \right|^\alpha - \left| \theta_1(x-1)^{d-1}\right|^\alpha - \left|\theta_2x^{d-1} \right|^\alpha dx.
\end{align}

\noindent The last above integral in expression \eqref{eq4_33} is the constant $C$ in (\ref{MRII_eq_6}) (see the proof of Lemma \ref{MRII_lema_3} for the value of $L$). Thus, $I(t) \to 0$, as $t \to \infty$. Therefore,
\begin{equation*} \label{MRII_eq_16}
	r(t) \sim -KI(t) \sim -KC|t|^{\alpha(d-1)+1}, \ \ \mbox{as} \ \ t \to \infty.
\end{equation*}
\vspace{-0.8cm}
\qed

\noindent The following theorem proves the long-range dependence property for the FIMA
Stable process.

\begin{teo} \label{MRII_teo_3}
	Let $\{Y_d(t)\}_{t \in \R}$ be the FIMA Stable process defined in (\ref{MR2_eq_1}). If $1 < \alpha < 2$ and $0 < d < 1 - \frac{1}{\alpha}$, then it has the property of long-range dependence or long-memory, in the sense of Definition \ref{MRII_def_1}.
\end{teo}

\noindent \textbf{Proof}
	Since $1<\alpha<2$ and $-1<d-1<-\frac{1}{\alpha}$ imply that $-2<\alpha(d-1) < -1$, we have by the Theorem \ref{MRII_teo_1} that $r(t)$ satisfies
	\begin{equation}
	|r(t)| = |r(\theta_1,\theta_2;t)| \sim KCt^{\alpha(d-1)+1} \; ,
	\end{equation}
	with $-1 < \alpha(d-1) + 1 < 0$.  \vspace{-0.2cm}\qed
	
\noindent Now, as an application of the representation Theorem \ref{MRII_teo_2} and the convergence Theorem \ref{MRI_teo_3}, we will prove a result which is basically a LLN
result for our setting. Suppose $\{t_j\}_{j \in \N}$ is a sequence of times such that $t_j \in \R$ and for every $j \in \N$ we consider $Y_d(t_j)$ with representation
\begin{equation}
	Y_d(t_j) = \int_{-\infty}^{t_j} g(t_j - x) dM_d(x) = \int_{\R} \I_{\{x \le t_j\}} g(t_j - x) dM_d(x) =: \int_{\R} g_j(x) dM_d(x) \; .
\end{equation}
\noindent Define the sequence of partial sums $\{S_n\}_{n \in \N}$ by
\begin{equation}\label{eq4_49}
	S_n := \sum_{j = 1}^{n} Y_d(t_j) = \sum_{j = 1}^{n} I_{M_d}(g_j) = I_{M_d}\left(\sum_{j = 1}^{n} g_j\right) = \int_{\R} \left[\sum_{j = 1}^{n} g_j(x)\right] dM_d(x)
\end{equation}
\noindent from the linearity of the integral $I_{M_d}$. So we state the following theorem.

\begin{teo}\label{TheoMRII4.5}
	If the sequence of times $\{t_j\}_{j \in \N}$ is a sequence of natural times, $t_j = j$, then
	\begin{equation} \label{MRII_eq_13}
		\dfrac{S_n}{n} \to 0 \text{ , as } n \to \infty \text{ in the } L^p(\Omega, \mathbb{P}) \text{ sense.}
	\end{equation}
\end{teo}

\noindent \textbf{Proof.} We start with some definitions. Let
	\begin{equation*}
		G_n(x) := \sum_{j = 1}^{n} g_j(x).
	\end{equation*}
	\noindent So the sequence in \eqref{eq4_49} becomes
	\begin{equation*}
		\dfrac{S_n}{n} = \int_{\R} \left[ \dfrac{G_n(x)}{n} \right] dM_d(x) \; .
	\end{equation*}
	Thus by Theorem \ref{MRI_teo_3} it suffices to show that $||G_n(x)/n||_1 \to 0$ and $||G_n(x)/n||_\alpha \to 0$, as $n \to \infty$, to conclude that
	\begin{equation*}
		\lim\limits_{n \to \infty} \int_{\R} \left[ \dfrac{G_n(x)}{n} \right] dM_d(x) = \int_{\R} 0 \; dM_d(x) = 0 \; .
	\end{equation*}
	We first note that as $g$ satisfies $|g(x)| \le C e^{-cx}$, for $x \ge 0$ and some constants $C,c > 0$,
	\begin{equation}\label{MRII_eq_14}
		\left|\frac{G_n(x)}{n}\right| \le \frac{1}{n}\sum_{j = 1}^{n} \I_{\{x \le j\}} |g(j - x)| \le \frac{C}{n}\sum_{j = 1}^{n} \I_{\{x \le j\}} e^{-c(j - x)}.
	\end{equation}
\noindent For every fixed $x \in \R$, the series on the right-hand side of the last equality in equation (\ref{MRII_eq_14}) converges to some positive constant limit $D_x$,
	\begin{equation*}
		\sum_{j = 1}^{\infty} \I_{\{x \le j\}} e^{-c(j - x)} = D_x.
	\end{equation*}
\noindent We get that for all $x \in \R$,
	\begin{equation*}
		\left|\frac{G_n(x)}{n}\right| \le \frac{CD_x}{n} \to 0 \text{ as } n \to \infty.
	\end{equation*}
	\noindent Furthermore, since $s \le t \Rightarrow \I_{\{x \le t\}} e^{-c(t - x)} \le \I_{\{x \le s\}} e^{-c(s - x)} \Rightarrow \sup_{j} \I_{\{x \le j\}} e^{-c(j - x)} = \I_{\{x \le 1\}} e^{-c(1 - x)}$ then, for every fixed $n \in \N$, by (\ref{MRII_eq_14}) we get that
	\begin{equation}
		\left|\frac{G_n(x)}{n}\right| \le \frac{C}{n}\sum_{j = 1}^{n} \I_{\{x \le 1\}} e^{-c(1 - x)} = C\I_{\{x \le 1\}} e^{-c(1 - x)} \; ,
	\end{equation}
	\noindent which clearly belongs to $L^1(\R) \cap L^\alpha(\R)$. We obtain the desired result by the Dominated Convergence theorem.
\vspace{-0.2cm}
\qed

\noindent The following corollary is a simple extension of the last result.

\begin{corol}
	If the sequence of times $\{t_j\}_{j \in \N}$ is bounded below $(t_j \in [a, \infty))$ such that $\lim\limits_{j \to \infty} j^{-\beta} \; t_j \ge K$ for some constants $\beta,K > 0$ then
	\begin{equation} \label{MRII_eq_15}
		\dfrac{S_n}{n} \to 0 \text{ , as } n \to \infty, \ \ \text{ in the } L^p(\Omega, \mathbb{P}) \text{ sense.}
	\end{equation}
\end{corol}

\noindent \textbf{Proof.} It suffices to note that since $\sup_{j} \I_{\{x \le t_j\}} e^{-c(t_j - x)} = \I_{\{x \le a\}} e^{-c(a - x)}$, for every $x \in \R$,
	\begin{equation*}
		\left|\frac{G_n(x)}{n}\right| \le \I_{\{x \le a\}} e^{-c(a - x)},
	\end{equation*}
	\noindent which also clearly belongs to $L^1(\R) \cap L^\alpha(\R)$. Since
$\lim\limits_{j \to \infty} j^{-\beta} \; t_j \ge K$, we get that for all $x \in \R$,
	\begin{equation*}
		\left|\frac{G_n(x)}{n}\right| \le \frac{C}{n}\sum_{j = 1}^{n} \I_{\{x \le Kj^\beta\}} e^{-c(Kj^\beta - x)} \le \frac{Ce^{cx}}{n}\sum_{j = 1}^{n} e^{-cKj^\beta} \le \frac{Ce^{cx}\tilde{D}}{n} \to 0 \text{ as } n \to \infty,
	\end{equation*}
	\noindent for
	\begin{equation*}
		\tilde{D} := \sum_{j = 1}^{\infty} e^{-cKj^\beta}.
	\end{equation*}
\vspace{-0.8cm}
\qed

\section{Conclusions}\label{conclusionsection}

We have successfully defined a continuous-time fractionally integrated moving average $\alpha$-stable process with the long-range dependence property. The results presented here can be applied in models that capture both high variability behavior and long-range dependence structure. We have shown an isometry allowing us to define stochastic integrals concerning the linear fractional stable motion by using the Riemann-Liouville fractional integrals. From this, a by-product was obtained as integration by parts formula. Theorem \ref{MRII_teo_2} reveals a core relationship between a FIMA Stable process and the LFSM, while Theorem \ref{MRII_teo_3} showed the long-range dependence property for the FIMA Stable processes. As an application of the isometry and integration by parts formula, law of large number's results follow immediately.

\vspace{0.3cm}

\subsection*{Acknowledgments}

G.L. Feltes was supported by CAPES-Brazil.
S.R.C. Lopes' research was partially supported by CNPq-Brazil.

\vspace{0.3cm}


\end{document}